\documentclass[12pt]{amsart}

\usepackage{amsmath}
\usepackage{amssymb}
\usepackage{amsfonts}
\usepackage{amsthm}
\usepackage{enumerate}
\usepackage{hyperref}
\usepackage{color}

\theoremstyle{plain}
\newtheorem{thm}{Theorem}[section]

\newtheorem{prop}[thm]{Proposition}

\newtheorem{cor}[thm]{Corollary}

\newtheorem{ques}[thm]{Question}

\newtheorem{conj}[thm]{Conjecture}

\theoremstyle{definition}
\newtheorem{dfn}[thm]{Definition}

\newtheorem{exmp}[thm]{Example}

\newtheorem{dfns-rems}[thm]{Definitions and Remarks}
\newtheorem{notas-rems}[thm]{Notations and Remarks}
\newtheorem{exmps-rems}[thm]{Examples and Remarks}

\begin{document}


\title[Stanley depth of powers of monomial ideals]{On the Stanley depth of powers of monomial ideals}


\author[S. A. Seyed Fakhari]{S. A. Seyed Fakhari}

\address{S. A. Seyed Fakhari, School of Mathematics, Statistics and Computer Science,
College of Science, University of Tehran, Tehran, Iran, and Institute of Mathematics, Vietnam Academy of Science and Technology, 18 Hoang Quoc Viet, Hanoi, Vietnam.}

\email{aminfakhari@ut.ac.ir}


\begin{abstract}
Let $\mathbb{K}$ be a field and $S=\mathbb{K}[x_1,\dots,x_n]$ be the
polynomial ring in $n$ variables over $\mathbb{K}$. In 1982, R. Stanley associated a combinatorial invariant to any finitely generated $\mathbb{Z}^n$-graded $S$-module which is now called Stanley depth. Stanley conjectured that this invariant is an upper bound for the depth of module. Stanley's conjecture has been disproved by Duval et al. \cite{abcj}, and the counterexample is a quotient of squarefree monomial ideals. On the other hand, there are evidences showing that Stanley's inequality can be true for high powers of monomial ideals. In this survey article, we collect the recent results in this direction. More precisely, we investigate the Stanley depth of powers, integral closure of powers and symbolic powers of monomial ideals.
\end{abstract}






\thanks{This research is partially funded by the Simons Foundation Grant Targeted for Institute of Mathematics, Vietnam Academy of Science and Technology.}


\maketitle


\section{Introduction} \label{sec1}

Let $\mathbb{K}$ be a field and let $S=\mathbb{K}[x_1,\dots,x_n]$
be the polynomial ring in $n$ variables over $\mathbb{K}$. Let
$M$ be a finitely generated $\mathbb{Z}^n$-graded $S$-module. Also, let
$u\in M$ be a homogeneous element and $Z\subseteq
\{x_1,\dots,x_n\}$. The $\mathbb {K}$-subspace $u\mathbb{K}[Z]$
generated by all elements $uv$, with $v$ a monomial in $\mathbb{K}[Z]$, is
called a {\it Stanley space} of dimension $|Z|$, if it is a free
$\mathbb{K}[Z]$-module. Here, as usual, $|Z|$ denotes the number
of elements of $Z$. A decomposition $\mathcal{D}$ of $M$ as a
finite direct sum of Stanley spaces is called a {\it Stanley
decomposition} of $M$. The minimum dimension of a Stanley space
in $\mathcal{D}$ is called the {\it Stanley depth} of
$\mathcal{D}$ and is denoted by ${\rm sdepth} (\mathcal {D})$.
The quantity $${\rm sdepth}(M):=\max\big\{{\rm sdepth}
(\mathcal{D})\mid \mathcal{D}\ {\rm is\ a\ Stanley\
decomposition\ of}\ M\big\}$$ is called the {\it Stanley depth}
of $M$. As a convention, we set ${\rm sdepth}(M)=\infty$, when $M$ is the zero module. For a reader friendly introduction to Stanley depth, we refer to
\cite{psty}.

\begin{exmp}
Consider the ideal $I=\langle x_1{x_2^2}, {x_1^2}x_2 \rangle$ in the polynomial ring $S=\mathbb{K}[x_1,x_2]$. Then $$\mathcal{D}_1:\ I=x_1x_2^2\mathbb{K}[x_2] \oplus x_1^2x_2\mathbb{K}[x_1,x_2]$$is a Stanley decomposition for $I$, with ${\rm sdepth}(\mathcal{D}_1)=1$. One can also write other Stanley decompositions for $I$. For example,$$\mathcal{D}_2:\ I=x_1^2x_2\mathbb{K} \oplus x_1^3x_2 \mathbb{K}[x_1]
\oplus x_1x_2^2 \mathbb{K}[x_2] \oplus {x_1^2}{x_2^2}\mathbb{K}[x_1,x_2],$$and$$\mathcal{D}_3:\ I=x_1^2x_2 \mathbb{K} \oplus x_1^3x_2 \mathbb{K} \oplus x_1^4x_2
\mathbb{K}[x_1] \oplus x_1x_2^2\mathbb{K}[x_2] \oplus {x_1^2}{x_2^2}
\mathbb{K}[x_1,x_2].$$It is clear that ${\rm sdepth}(\mathcal{D}_2)={\rm sdepth}(\mathcal{D}_3)=0$. It follows from the definition of Stanley depth that ${\rm sdepth}(I)\geq 1$ and it can be easily verified that indeed the equality holds, i.e., ${\rm sdepth}(I)=1$.
\end{exmp}

We say that a $\mathbb{Z}^n$-graded $S$-module $M$ satisfies {\it Stanley's inequality} if $${\rm depth}(M) \leq
{\rm sdepth}(M).$$ In fact, Stanley \cite{s} conjectured that the above inequality holds for every finitely generated $\mathbb{Z}^n$-graded $S$-module. Stanley's conjecture has been disproved by Duval, Goeckner, Klivans and Martin \cite{abcj}. In fact they construct a non-partitionable Cohen-Macaulay simplicial complex, and then using a result of Herzog, Soleyman Jahan and Yassemi \cite[Corollary 4.5]{hsy} deduce that the Stanley Reisner ring of this simplicial complex does not satisfy Stanley's inequality. In particular, the counterexample given in \cite{abcj} lives in the category of squarefree monomial ideals. Thus, one can still ask whether Stanley's inequality holds for non-squarefree monomial ideals. Of particular interest is the validity of Stanley's inequality for high powers of monomial ideals. In this survey article we review the recent developments in this regard. In 2013, Herzog \cite{h} published his nice survey on Stanley depth. In fact we complement his survey by collecting the results obtained since then with focus on powers of monomial ideals.


\section{Ordinary powers} \label{sec2}

In this section, we consider the ordinary powers of monomial ideal. As we explained in introduction, it is natural to ask whether the high powers of monomial ideals satisfy Stanley's inequality. In fact, this question was posed in \cite{s5}.

\begin{ques} [\cite{s5}, Question 1.1] \label{q1}
Let $I$ be a monomial ideal. Is it true that $I^k$ and $S/I^k$ satisfy Stanley's inequality for every integer $k\gg 0$?
\end{ques}

In the following subsections we will see that Question \ref{q1} has positive answer when $I$ belongs to interesting classes of monomial ideals.

\subsection{Maximal ideal and complete intersections}

Let $\frak{m}=(x_1,\ldots,x_n)$ denote the maximal ideal of $S$. It is clear that for every integer $k\geq 1$, ${\rm depth}(S/\frak{m}^k)=0$. Hence, $S/\frak{m}^k$ satisfies Stanley's inequality, for any $k\geq 1$. Indeed, since $S/\frak{m}^k$ is an Artinian ring, we also have ${\rm sdepth}(S/\frak{m}^k)=0$, for every integer $k\geq 1$. On the other hand, ${\rm depth}(\frak{m}^k)=1$ and by \cite[Corollary 24]{h}, we know that the Stanley depth of any monomial ideal is at least one. This implies that $\frak{m}^k$ satisfies Stanley's inequality for every integer $k\geq 1$. However, computing the exact value of Stanley depth of $\frak{m}^k$ is not easy. In 2010, Bir\'o, Howard, Keller, Trotter and Young \cite{bhkty} proved that ${\rm sdepth}(\frak{m})=\lceil n/2\rceil$. Cimpoea\c{s} \cite{c} determined an upper bound for the Stanley depth of powers of $\frak{m}$. More precisely, he proved the following results.

\begin{thm} [\cite{c}, Theorem 2.2] \label{maximal}
For every integer $k\geq 1$, we have$${\rm sdepth}(\frak{m}^k)\leq \bigg\lceil \frac{n}{k+1}\bigg\rceil.$$In particular, for every integer $k\geq n-1$, we have ${\rm sdepth}(\frak{m}^k)=1$.
\end{thm}

Cimpoea\c{s} \cite{c} also conjectured that the inequality obtained in the above theorem is indeed equality, i.e.,$${\rm sdepth}(\frak{m}^k)=\bigg\lceil \frac{n}{k+1}\bigg\rceil,$$for every $k\geq 1$.

In 2018, Cimpoea\c{s} \cite{c1} extended Theorem \ref{maximal} by determining bounds for the Stanley depth of complete intersection monomial ideals.

\begin{thm} [\cite{c1}, Proposition 2.14 and Theorem 2.15] \label{compint}
Let $I$ be a complete intersection monomial ideal which is minimally generated by $t$ monomials.

\begin{itemize}
\item[(i)] For every integer $k\geq 1$, we have$$n-t+1\leq {\rm sdepth}(I^k)\leq n-t+\bigg\lceil \frac{t}{k+1}\bigg\rceil.$$In particular, if $k\geq t-1$, then ${\rm sdepth}(I^k)=n-t+1$.
\item[(ii)] For every integer $k\geq 1$, we have$${\rm sdepth}(S/I^k)={\rm sdepth}(I^k/I^{k+1})=\dim(S/I)=n-t.$$
\end{itemize}
\end{thm}

As an immediate consequence of Theorem \ref{compint}, we conclude that for any complete intersection monomial ideal and every integer $k\geq 1$, the modules, $I^k, S/I^k$ and $I^k/I^{k+1}$ satisfy Stanley's inequality. In particular, Question \ref{q1} has positive answer in this case.

\subsection{Polymatroidal ideals}

We begin this subsection by recalling the definition of polymatroidal ideals.
\begin{dfn}
Let $I$ be a monomial ideal of $S$ which is
generated in a single degree and assume that $G(I)$ is the set of minimal
monomial generators of $I$. The ideal $I$ is called {\it polymatroidal} if
the following exchange condition is satisfied: For monomials $u=x^{a_1}_1
\ldots x^{a_n}_ n$ and $v = x^{b_1}_1\ldots x^{b_n}_ n$ belonging to $G(I)$
and for every $i$ with $a_i > b_i$, one has $j$ with $a_j < b_j$ such that
$x_j(u/x_i)\in G(I)$.
\end{dfn}

Weakly polymatroidal ideals are generalization of polymatroidal ideals and
they are defined as follows.

\begin{dfn} [\cite{mm}, Definition 1.1]
A monomial ideal $I$ of $S$ is called {\it weakly
polymatroidal} if for every two monomials $u = x_1^{a_1} \ldots x_n^{a_n}$
and $v = x_1^{b_1} \ldots x_n^{b_n}$ in $G(I)$ such that $a_1 = b_1,
\ldots, a_{t-1} = b_{t-1}$ and $a_t > b_t$ for some $t$, there exists $j >
t$ such that $x_t(v/x_j)\in I$.
\end{dfn}

It is obvious that any polymatroidal ideal is weakly polymatroidal.

Let $I$ be a weakly polymatroidal ideal. In \cite[Theorem 2.4]{s7}, we proved that $S/I$ satisfies Stanley's inequality. We also know from
\cite[Theorem 12.6.3]{hh'} that every power of a polymatroidal ideal is again a polymatroidal ideal. As a consequence, for any polymatroidal ideal $I$ and any integer $k\geq 1$, the module $S/I^k$ satisfies Stanley's inequality. It is natural to ask whether $I^k$ satisfies Stanley's inequality. Before answering this question, we recall the concept of having linear quotiens, introduced in \cite{ht1}.

\begin{dfn}
Let $I$ be a monomial ideal and assume that $G(I)$ is the set of minimal
monomial generators of $I$. We say that $I$ has linear quotients if there is a linear order $u_1\prec u_2 \prec \ldots \prec u_m$ on $G(I)$, with the property that for every $2\leq i\leq m$, the ideal $(u_1, \ldots, u_{i-1}):u_i$ is generated by a subset of the variables.
\end{dfn}

Soleyman Jahan \cite{so} proves that Stanley's inequality holds for any monomial ideal which has linear quotients. On the other hand, by \cite[Theorem 1.3]{mm}, we know that any weakly polymatroidal ideal has linear quotients. This implies that every weakly polymatroidal ideal satisfies Stanley's inequality. Since every power of a polymatroidal ideal is again a polymatroidal ideal, we deduce that for any polymatroidal ideal $I$ and any integer $k\geq 1$, the ideal $I^k$ satisfies Stanley's inequality.

By the above argument, we know that Question \ref{q1} has positive answer for polymatroidal ideals. This result was also obtained in \cite{psy1}.

Let $I$ be a monomial ideal of $S$ with Rees algebra $\mathcal{R}(I)=\bigoplus_
{k=0}^{\infty}I^k$. The $\mathbb{K}$-algebra $\mathcal{R}(I)/\mathfrak{m}\mathcal{R}(I)$ is
called the {\it fibre ring} and its Krull dimension is called the {\it
analytic spread} of $I$, denoted by $\ell(I)$. A classical result by Burch \cite{b'} states
that $$\min_k{\rm depth}(S/I^k)\leq n-\ell(I).$$By a theorem of
Brodmann \cite{b}, ${\rm depth}(S/I^k)$ is constant for large $k$. We call
this constant value the {\it limit depth} of $I$, and denote it by
$\lim_{k\rightarrow \infty}{\rm depth}(S/I^k)$. Brodmann improved the Burch's
inequality by showing that$$\lim_{k\rightarrow \infty}{\rm depth}(S/I^k)
\leq n-\ell(I).$$

We know from \cite[Corollary 3.5]{hrv} that equality occurs in the above inequality, if $I$ is a polymatroidal ideal. In fact, we will see in the next section that equality holds in Burch's inequality for a larger class of ideals, namely, the class of {\it normal} ideals.

Inspired by the limit behavior of depth of powers of ideals, Herzog \cite{h} proposed the following conjecture.

\begin{conj} [\cite{h}, Conjecture 59] \label{limsd}
For every monomial ideal $I$, the sequence $\{{\rm sdepth}(I^k)\}_{k=1}^{\infty}$ is convergent.
\end{conj}

This conjecture is widely open. However, by Theorem \ref{compint}, it has positive answer for complete intersections. Also, we will see in Section \ref{sec4} that the assertion of this conjecture is true for any normally torsionfree squarefree monomial ideal.

Let $I$ be a weakly polymatroidal ideal which is generated in a single degree. We know from \cite[Theorem 2.5]{psy1} that ${\rm depth}(S/I)\geq n-\ell(I)$. Since $I$ and $S/I$ satisfy Stanley's inequality, it follows that$${\rm sdepth}(S/I)\geq n-\ell(I) \ \ \ \ \ {\rm and} \ \ \ \ \ {\rm sdepth}(I)\geq n-\ell(I)+1.$$Restricting to the class of polymatroidal ideals, for any integer $k\geq 1$ and any polymatroidal ideal $I$, we have$${\rm sdepth}(S/I^k)\geq n-\ell(I) \ \ \ \ \ {\rm and} \ \ \ \ \ {\rm sdepth}(I^k)\geq n-\ell(I)+1.$$Indeed, we expect that the equality holds in the above inequality for every $k\gg 0$. In other words, not only we believe that Conjecture \ref{limsd} is true for every polymatroidal ideal $I$, but we also have a prediction for the limit value of the Stanley depth of powers of $I$.

\begin{conj}
Let $I$ be a polymatroidal ideal. Then$${\rm sdepth}(S/I^k)=n-\ell(I) \ \ \ \ \ {\rm and} \ \ \ \ \ {\rm sdepth}(I^k)=n-\ell(I)+1,$$for any integer $k\gg 0$.
\end{conj}

\subsection{Edge ideals}

There is a natural correspondence between quadratic squarefree monomial ideals of $S$ and finite simple graphs with $n$ vertices. To every simple graph $G$ with vertex set $V(G)=\big\{x_1, \ldots, x_n\big\}$ and edge set $E(G)$, we associate its {\it edge ideal} $I=I(G)$ defined by
$$I(G)=\big(x_ix_j: x_ix_j\in E(G)\big)\subseteq S.$$

Stanley depth of powers of edge ideals has been studied in \cite{asy}, \cite{fm}, \cite{psy} and \cite{s3}. Before reviewing the main results of these papers, we mention the following result of Trung, concerning the depth of high powers of edge ideals.

\begin{thm} [\cite{t}, Theorems 4.4 and 4.6] \label{trung}
Let $G$ be a graph with $n$ vertices and $p$ bipartite connected components. Then for every integer $k\geq n-1$, we have$${\rm depth}(S/I(G)^k)=p.$$
\end{thm}

Note that by \cite[Page 50]{v}, for every graph $G$ with $n$ vertices and $p$ bipartite connected components, we have $\ell(I(G))=n-p$. Thus, Theorem \ref{trung}, essentially say that$$\lim_{k\rightarrow \infty}{\rm depth}(S/I(G)^k)
=n-\ell(I(G)),$$i.e., equality occurs in Burch's inequality.

Pournaki, Yassemi and the author \cite{psy} studied the Stanley depth of $S/I(G)^k$, where $G$ is a forest (i.e., a graph with no cycle). They proved that for every forest with $p$ connected components and any integer $k\geq 1$, we have$${\rm sdepth}(S/I(G)^k)\geq p.$$ This together with Theorem \ref{trung} implies that for any forest $G$ with $n$ vertices, the module $S/I(G)^k$ satisfies Stanley's inequality for any integer $k\geq n-1$. This result was then extended in \cite{s3}, to any arbitrary graph, as follows.

\begin{thm} [\cite{s3}, Theorem 2.3 and Corollary 2.5] \label{quo}
Let $G$ be a graph with $n$ vertices and $p$ bipartite connected components. Then for every integer $k\geq 1$, we have ${\rm sdepth}(S/I(G)^k)\geq p$. In particular, $S/I(G)^k$ satisfies Stanley's inequality for any integer $k\geq n-1$.
\end{thm}

We know from the above theorem that for any graph $G$, the module $S/I(G)^k$ satisfies Stanley's inequality for $k\gg 0$. But how about $I(G)^k$? By Theorem \ref{trung}, in order to prove Stanley's inequality for high powers of $I(G)$, we need to prove ${\rm sdepth}(I(G)^k)\geq p+1$, for every integer $k\gg 0$. We do know whether this inequality holds for any arbitrary graph. However, we have a partial result, as follows. We recall that for any graph $G$ and every subset $U\subset V(G)$, the graph $G\setminus U$ has vertex set $V(G\setminus U)=V(G)\setminus U$ and edge set $E(G\setminus U)=\{e\in E(G)\mid e\cap U=\emptyset\}$.

\begin{thm} [\cite{s3}, Theorem 3.1] \label{ideal}
Let $G$ be a graph and assume that $H$ is a connected component of $G$ with at least one edge. Suppose that $h$ is the number of bipartite connected components of $G\setminus V(H)$. Then for every integer $k\geq 1$, we have$${\rm sdepth}(I(G)^k)\geq \min_{1\leq l \leq k}\{{\rm sdepth}_{S'}(I(H)^l)\}+h,$$where $S'=\mathbb{K}[x_i\mid x_i\in V(H)]$.
\end{thm}

Assume that $G$ has a non-bipartite connected component and call it $H$. Then by \cite[Corollary 24]{h}, for every integer $l\geq 1$, we have ${\rm sdepth}(I(H)^l)\geq 1$. Thus, it follows from Theorem \ref{ideal} that in this case, ${\rm sdepth}(I(G)^k)\geq p+1$, where $p$ is the number of bipartite connected components of $G$ and $k\geq 1$ is an arbitrary positive integer. Assume now that $G$ is a bipartite graph. Using Theorem \ref{ideal}, in order to prove the inequality ${\rm sdepth}(I(G)^k)\geq p+1$, it is enough to prove it only for the class of connected bipartite graphs. Thus, we raise the following question.

\begin{ques} [\cite{s3}, Question 3.3] \label{questbip}
Let $G$ be a connected bipartite graph (with at least one edge) and suppose $k\geq 1$ is an integer. Is it true that ${\rm sdepth}(I(G)^k)\geq 2$?
\end{ques}

We investigated this question in \cite{s8} and proved that it has positive answer for small $k$. More precisely, we proved the following result.

\begin{thm} [\cite{s8}, Theorem 3.4] \label{girth}
Let $G$ be a connected bipartite graph (with at least one edge) and let $g$ be a positive integer. Suppose $G$ has no cycle of length at most $g-1$. Then for every positive integer $k\leq g/2+1$, we have ${\rm sdepth}(I(G)^k)\geq 2$.
\end{thm}

Theorem \ref{girth}, in particular implies that ${\rm sdepth}(I(G)^k)\geq 2$, for any integer $k\geq 1$, provided that $G$ is a tree (i.e., a connected forest). Combining this result with Theorem \ref{ideal} implies that if $G$ is a bipartite graph and at least one of the connected components of $G$ is a tree, then for every integer $k\geq 1$, we have ${\rm sdepth}(I(G)^k)\geq p+1$, where $p$ is the number of (bipartite) connected components of $G$. All of all, we obtain the following theorem.

\begin{thm} [\cite{s3}, Corollary 3.6] \label{sin}
Assume that $G$ is a graph with $n$ vertices, such that
\begin{itemize}
\item[(i)] $G$ is a non-bipartite graph, or
\item[(ii)] at least one of the connected components of $G$ is a tree with at least one edge.
\end{itemize}
Then for every integer $k\geq n-1$, the ideal $I(G)^k$ satisfies Stanley's inequality.
\end{thm}

Let $I$ be a monomial ideal. We know by \cite[Theorem 1.2]{hh''} that the sequence $\{{\rm depth}(I^k/I^{k+1})\}_{k=1}^{\infty}$ is convergent and moreover,$$\lim_{k\rightarrow\infty}{\rm depth}(I^k/I^{k+1})= \lim_{k\rightarrow\infty}{\rm depth}(S/I^k).$$Therefore, using Theorem \ref{trung}, we conclude that for any graph $G$,$$\lim_{k\rightarrow\infty}{\rm depth}(I(G)^k/I(G)^{k+1})=p,$$where $p$ is the number of bipartite connected components of $G$. In \cite{s3}, we also studied the Stanley depth of $I(G)^k/I(G)^{k+1}$ and proved that it satisfies Stanley's inequality for any $k\gg 0$. In fact, we proved the following result.

\begin{thm} [\cite{s3}, Theorem 2.2 and Corollary 2.6] \label{twoquo}
Let $G$ be a graph and suppose $p$ is the number of bipartite connected components of $G$. Then for every integer $k\geq 0$, we have ${\rm sdepth}(I(G)^k/I(G)^{k+1})\geq p$. In particular, $I(G)^k/I(G)^{k+1}$ satisfies Stanley's inequality, for every integer $k\gg 0$.
\end{thm}

We mention that in the special case, when $G$ is a forest, Theorem \ref{twoquo} was proved in \cite[Theorem 3.1]{asy}.

The {\it diameter} of a connected graph is the maximum distance between any two
vertices. Here, the {\it distance} between two vertices is the minimum
length of a path connecting the vertices.

Fouli and Morey \cite{fm} studied the Stanley depth of small powers of edge ideals and determined a lower bound for it.

\begin{thm} [\cite{fm}, Theorem 4.18] \label{diam}
Assume that $G$ is a graph with $c$ connected components and let $d$ denote the maximum of the diameters of the connected components of $G$. Then for every integer $1\leq t\leq 3$, we have$${\rm sdepth}(S/I(G)^t)\geq \bigg\lceil\frac{d-4t+5}{3}\bigg\rceil+c-1.$$
\end{thm}

Fouli and Morey \cite[Corollary 3.3, Theorems 4.4 and 4.13]{fm} also show that the inequality of Theorem \ref{diam} remains true, if one replaces sdepth with depth.


\section{Integral closure of powers} \label{sec3}

The study of Stanley depth of integral closure of powers of monomial ideals was initiated in \cite{s1} and continued in \cite{s8}. Before stating the results of theses paper, we recall some definitions and basic facts from the theory of integral closure.

Let $I\subset S$ be an arbitrary ideal. An element $f \in S$ is
{\it integral} over $I$, if there exists an equation
$$f^k + c_1f^{k-1}+ \ldots + c_{k-1}f + c_k = 0 {\rm \ \ \ \ with} \ c_i\in I^i.$$
The set of elements $\overline{I}$ in $S$ which are integral over $I$ is the {\it integral closure}
of $I$. It is known that the integral closure of a monomial ideal $I\subset S$ is a monomial ideal
generated by all monomials $u \in S$ for which there exists an integer $k$ such that
$u^k\in I^k$ (see \cite[Theorem 1.4.2]{hh'}). The ideal $I$ is {\it integrally closed}, if $I = \overline{I}$, and $I$ is {\it normal} if all powers
of $I$ are integrally closed. By \cite[Theorem 3.3.18]{v'}, a monomial ideal $I$ is normal if and only if the Rees algebra $\mathcal{R}(I)$ is a normal ring.

We first notice that there is no general inequality between the Stanley depth of $S/I$ and that of $S/\overline{I}$. This will be illustrated in the following examples.

\begin{exmp} [\cite{s1}, Example 1.2] \label{ex1}
Let $I=(x_1^2, x_2^2, x_1x_2x_3)$ be a monomial ideal in the polynomial ring $S=\mathbb{K}[x_1, x_2,x_3]$. It is not difficult to see that $\overline{I}=(x_1^2, x_2^2,x_1x_2)$. Then the maximal ideal $\mathfrak{m}=(x_1,x_2,x_3)$ is an associated prime of $S/I$ and it follows from \cite[Proposition 1.3]{hvz} that ${\rm sdepth}(S/I)=0$. Since $\mathfrak{m}$ is not an associated prime of $S/\overline{I}$, it follows from \cite[Proposition 2.13]{bku} that ${\rm sdepth}(S/\overline{I})\geq 1$. Thus, in this example ${\rm sdepth}(S/I) < {\rm sdepth}(S/\overline{I})$.
\end{exmp}

\begin{exmp} [\cite{s1}, Example 1.3] \label{ex2}
Let $I=(x_1^2x_2^2, x_1^2x_3^2, x_2^2x_3^2)$ be a monomial ideal in the polynomial ring $S=\mathbb{K}[x_1, x_2, x_3]$. The maximal ideal $\mathfrak{m}=(x_1,x_2,x_3)$ is not an associated prime of $S/I$ and hence, using \cite[Proposition 2.13]{bku}, we conclude that ${\rm sdepth}(S/I)\geq 1$. On the other, we know from \cite[Theorem 2.4]{j} that $\mathfrak{m}$ is an associated prime of $S/\overline{I}$ and therefore, \cite[Proposition 1.3]{hvz} implies that ${\rm sdepth}(S/\overline{I})=0$. Thus, in this example ${\rm sdepth}(S/I) > {\rm sdepth}(S/\overline{I})$.
\end{exmp}

Although there is no general inequality between ${\rm sdepth}(S/I)$ and ${\rm sdepth}(S/\overline{I})$, but we will see in the following theorem that the Stanley depth of $S/\overline{I}$ provides an upper bound for the Stanley depth of the quotient ring of some powers of $I$.

\begin{thm} [\cite{s1}, Theorem 2.8] \label{int1}
Let $I_2\subseteq I_1$ be two monomial ideals in $S$. Then there  exists an integer $k\geq 1$, such that for every $s\geq 1$
$${\rm sdepth} (I_1^{sk}/I_2^{sk}) \leq {\rm sdepth} (\overline{I_1}/\overline{I_2}).$$
\end{thm}

In particular, we have the following corollary.

\begin{cor} \label{mainc}
Let $I\subset S$ be a monomial ideal. Then there exist integers $k_1, k_2 \geq 1$, such that for every $s\geq 1,$
$${\rm sdepth} (I^{sk_1}) \leq {\rm sdepth} (\overline{I})$$
and
$${\rm sdepth} (S/I^{sk_2}) \leq {\rm sdepth} (S/\overline{I}).$$
\end{cor}

We mention that the assertions of Corollary \ref{mainc} remain true if one replaces sdepth with depth, \cite[Theorem 4.5]{s8}.

In Question \ref{q1}, we asked whether the high powers of an ideal satisfy Stanley's inequality. One can ask a similar question by  replacing $I^k$ with its integral closure. This question is posed in \cite{s8}.

\begin{ques} [\cite{s8}, Question 1.2] \label{q2}
Let $I$ be a monomial ideal. Is it true that $\overline{I^k}$ and $S/\overline{I^k}$ satisfy Stanley's inequality for every integer $k\gg 0$?
\end{ques}

Before we focus on the above question, we recall the following result of Hoa and Trung concerning the depth of integral closure of high powers of monomial ideals.

\begin{thm} [\cite{ht}, Lemma 1.5] \label{depthint}
Let $I$ be a monomial ideal of $S$. Then ${\rm depth}(S/\overline{I^k})=n-\ell(I)$, for every integer $k\gg 0$.
\end{thm}

According to the above theorem, Question \ref{q2} is equivalent to the following question.

\begin{ques} [\cite{s8}, Question 1.3] \label{q3}
Let $I$ be a monomial ideal. Is it true that the inequalities ${\rm sdepth}(\overline{I^k})\geq n-\ell(I)+1$ and ${\rm sdepth}(S/\overline{I^k})\geq n-\ell(I)$ hold, for every integer $k\gg 0$?
\end{ques}

Let $I$ be a monomial ideal of $S$ and assume that ${\rm sdepth}(S/I^k)\geq n-\ell(I)$ (resp. ${\rm sdepth}(I^k)\geq n-\ell(I)+1$), for every integer $k\gg 0$. It follows from Corollary \ref{mainc} that ${\rm sdepth}(S/\overline{I^k})\geq n-\ell(I)$ (resp. ${\rm sdepth}(\overline{I^k})\geq n-\ell(I)+1$), for every integer $k\gg 0$. Thus, the answers of Questions \ref{q2} and \ref{q3} are positive for $I$. This argument together with Theorem \ref{compint} implies the following result, concerning the Stanley depth of integral closure of powers complete intersection monomial ideals.

\begin{thm}
Let $I$ be a complete intersection monomial ideal which is minimally generated by $t$ monomials.

\begin{itemize}
\item[(i)] For every integer $k\geq 1$, we have${\rm sdepth}(\overline{I^k})\geq n-t+1$.
\item[(ii)] For every integer $k\geq 1$, we have$${\rm sdepth}(S/\overline{I^k})=n-t.$$
\end{itemize}
\end{thm}

Note that in part (ii) of the above theorem, we use the fact that for any complete intersection monomial ideal and any integer $k\geq 1$, the dimension of $S/\overline{I^k}$ is $n-t$, where $t$ is the number of minimal monomial generators of $I$.

Restricting to edge ideals, combining the above argument with Theorems \ref{quo} and \ref{sin} implies the following results.

\begin{thm} [\cite{s8}, Theorem 3.2] \label{iquo}
Let $G$ be a graph and suppose that $p$ is the number of bipartite connected components of $G$. Then for every integer $k\geq 1$, we have ${\rm sdepth}(S/\overline{I(G)^k})\geq p$. In particular, $S/\overline{I(G)^k}$ satisfies Stanley's inequality for every integer $k\gg 0$.
\end{thm}

\begin{thm} [\cite{s8}, Theorem 3.3] \label{isin}
Let $G$ be a non-bipartite graph and suppose that $p$ is the number of bipartite connected components of $G$. Then for every integer $k\geq 1$, we have ${\rm sdepth}(\overline{I(G)^k})\geq p+1$. In particular, $\overline{I(G)^k}$ satisfies Stanley's inequality for every integer $k\gg 0$.
\end{thm}

Assume that $G$ is a bipartite graph. We know from \cite[Theorem 1.4.6 and Corollary 10.3.17]{hh'} that for any integer $k\geq 1$, the equality $I(G)^k=\overline{I(G)^k}$ holds. Therefore, $\overline{I(G)^k}$ satisfies Stanley's inequality if and only if $I(G)^k$ satisfies that inequality. Because of this reason, we exclude the case of bipartite graphs in Theorem \ref{isin}.

Let $I$ be a monomial ideal. It is also reasonable to study the depth and the Stanley depth of $\overline{I^k}/\overline{I^{k+1}}$. In \cite{s8}, we proved the following result about the depth of these modules for large $k$.

\begin{thm} [\cite{s8}, Theorem 4.1] \label{dtwoquo}
For any nonzero monomial ideal $I\varsubsetneq S$, the sequence $\{{\rm depth}(\overline{I^k}/\overline{I^{k+1}})\}_{k=0}^{\infty}$ is convergent and moreover,$$\lim_{k\rightarrow\infty}{\rm depth}(\overline{I^k}/\overline{I^{k+1}})=n-\ell(I).$$
\end{thm}

According to Theorem \ref{dtwoquo}, in order to prove that $\overline{I^k}/\overline{I^{k+1}}$ satisfies Stanley's inequality, for $k\gg 0$, we must show that ${\rm sdepth}(\overline{I^k}/\overline{I^{k+1}})\geq n-\ell(I)$, for high $k$.

Let $I$ be a monomial ideal of $S$ with ${\rm sdepth}(I^k/I^{k+1})\geq n-\ell(I)$, for every integer $k\gg 0$, say for $k\geq k_0$. We fix an integer $k\geq 1$. By Corollary \ref{mainc}, there exists an integer $s$ with $sk\geq k_0$ such that$${\rm sdepth}(\overline{I^k}/\overline{I^{k+1}})\geq {\rm sdepth}(I^{sk}/I^{s(k+1)}).$$On the other hand, as $\mathbb{K}$-vector spaces, we have$$I^{sk}/I^{s(k+1)}=\bigoplus_{i=sk}^{sk+s-1}I^i/I^{i+1}.$$By the definition of Stanley depth we conclude that$${\rm sdepth}(I^{sk}/I^{s(k+1)})\geq \min\big\{{\rm sdepth}(I^i/I^{i+1}) \mid i=sk, \ldots, sk+s-1\big\}\geq n-\ell(I),$$ where the last inequality follows from the assumption. Therefore, $${\rm sdepth}(\overline{I^k}/\overline{I^{k+1}})\geq n-\ell(I).$$Hence, $\overline{I^k}/\overline{I^{k+1}}$ satisfies Stanley's inequality, for $k\gg 0$. In particular cases, it follows from Theorems \ref{compint} and \ref{twoquo} that $\overline{I^k}/\overline{I^{k+1}}$ satisfies Stanley's inequality, for every integer $k\gg 0$, if $I$ is either a complete intersection monomial ideal or an edge ideal.

Let $I$ be a normal ideal. By \cite[Proposition 10.3.2]{hh'},$$\lim_{k\rightarrow \infty}{\rm depth}(S/I^k)=n-\ell(I).$$ Hence, if $I^k$ and $S/I^k$ satisfy Stanley's inequality for large $k$, we must have$${\rm sdepth}(S/I^k)\geq n-\ell(I) \ \ \ {\rm and} \ \ \ {\rm sdepth}(I^k)\geq n-\ell(I)+1.$$In fact, in \cite{s1}, we conjectured that the above inequalities hold in a more general setting.

\begin{conj} [\cite{s1}, Conjecture 2.6] \label{conje}
Let $I\subset S$ be an integrally closed monomial ideal. Then ${\rm sdepth}(S/I)\geq n-\ell(I)$ and ${\rm sdepth} (I)\geq n-\ell(I)+1$.
\end{conj}

The following example shows that the inequalities of Conjecture \ref{conje} do not necessarily hold, if $I$ is not integrally closed.

\begin{exmp} [\cite{s1}, Example 2.5] \label{example}
Consider the ideal $I=(x_1^2, x_2^2, x_1x_2x_3, x_1x_2x_4)$ in the polynomial ring $S=\mathbb{K}[x_1, x_2, x_3, x_4]$. Then $\ell(I)=2$. But $\mathfrak{m}=(x_1, x_2,x_3, x_4)$ is an associated prime of $S/I$ and therefore, we conclude from \cite[Proposition 1.3]{hvz} that ${\rm sdepth}(S/I)=0$ and by \cite[Corollary 1.2]{i}, ${\rm sdepth}(I)\leq 2$. This shows that the inequalities ${\rm sdepth}(S/I)\geq n-\ell(I)$ and ${\rm sdepth} (I)\geq n-\ell(I)+1$ do not hold for $I$.
\end{exmp}

As we mentioned in Section \ref{sec2}, the inequalities of Conjecture \ref{conje} are true for any polymatroidal ideal (we know from \cite[Theorem 3.4]{hrv} that any polymatroidal ideal is integrally closed). Also, in \cite[Corollary 3.4]{s2}, we verified Conjecture \ref{conje} for any squarefree monomial ideal which is generated in a single degree.

We close this section by the following result which permits us to compare the Stanley depth of integral closure of a monomial ideal and its powers.

\begin{thm} [\cite{s1}, Theorem 2.8] \label{ipow}
Let $J\subseteq I$ be two monomial ideals in $S$. Then for every integer $k\geq 1$
$${\rm sdepth} (\overline{I^k}/\overline{J^k}) \leq {\rm sdepth} (\overline{I}/\overline{J}).$$
\end{thm}

The following corollary is an immediate consequence of Theorem \ref{ipow}.

\begin{cor} \label{ipow1}
Let $I \subset S$ be a monomial ideal. Then for every integer $k\geq 1$,
$${\rm sdepth} (\overline{I^k}) \leq {\rm sdepth} (\overline{I})$$
and
$${\rm sdepth} (S/\overline{I^k}) \leq {\rm sdepth} (S/\overline{I}).$$
\end{cor}

We mention that the inequalities of Corollary \ref{ipow1} remain true, if one replaces sdepth with depth and this has been proved by Hoa and Trung \cite[Lemma 2.5]{ht}.


\section{Symbolic powers} \label{sec4}

In this section, we collect the recent results concerning the Stanley depth of symbolic powers of squarefree monomial ideals. We first recall the definition of symbolic powers and then we continue in two subsections.

\begin{dfn}
Let $I$ be an ideal of $S$ and let ${\rm Min}(I)$ denote the set of minimal primes of $I$. For every integer $k\geq 1$, the $k$-th {\it symbolic power} of $I$,
denoted by $I^{(k)}$, is defined to be$$I^{(k)}=\bigcap_{\frak{p}\in {\rm Min}(I)} {\rm Ker}(S\rightarrow (S/I^k)_{\frak{p}}).$$
\end{dfn}

Let $I$ be a squarefree monomial ideal in $S$ and suppose that $I$ has the irredundant
primary decomposition $$I=\frak{p}_1\cap\ldots\cap\frak{p}_r,$$ where each
$\frak{p}_i$ is a prime ideal generated by a subset of the variables of
$S$. It follows from \cite[Proposition 1.4.4]{hh'} that for every integer $k\geq 1$, $$I^{(k)}=\frak{p}_1^k\cap\ldots\cap
\frak{p}_r^k.$$

\subsection{Asymptotic behavior of Stanley depth of symbolic powers}

Let $I$ be a squarefree monomial ideal. As we mentioned in Section \ref{sec2}, based on the limit behavior of depth of powers of $I$, Herzog \cite{h} conjectured that the Stanley depth of $S/I^k$ is constant for large $k$ (see Conjecture \ref{limsd}). On the other hand, it is known that if one replaces the ordinary powers by symbolic powers, then again the depth function stabilizes. In fact, Hoa, Kimura, Terai and Trung \cite{hktt} are even able to compute the limit value of this function. In order to state their result, we need the following definition.

\begin{dfn}
Suppose $I$ is a squarefree monomial ideal and let $\mathcal{R}_s(I)=\bigoplus_
{k=0}^{\infty}I^{(k)}$ be the {\it symbolic Rees ring} of $I$. The Krull dimension of $\mathcal{R}_s(I)/{{\frak{m}}\mathcal{R}(I)}$ is called the {\it symbolic analytic spread} of $I$ and is denoted by $\ell_s(I)$.
\end{dfn}

Let $I$ be a squarefree monomial ideal. Varbaro \cite[Proposition 2.4]{v1} showed that$$\min_k{\rm depth}(S/I^{(k)})=n-\ell_s(I).$$ In \cite{hktt},  Hoa, Kimura, Terai and Trung proved that the minimum and the limit of the sequence $\{{\rm depth}(S/I^{(k)})\}_{k=1}^{\infty}$ coincide. Indeed, they showed the following stronger result. In the following theorem, ${\rm bight}(I)$ denotes the maximum height of associated primes of $I$.

\begin{thm} [\cite{hktt}, Theorem 2.4] \label{dsymlim}
Let $I$ be a squarefree monomial ideal of $S$. Then ${\rm depth}(S/I^{(k)})= n-\ell_s(I)$, for every integer $k\geq n(n+1){\rm bight}(I)^{n/2}$.
\end{thm}

As the depth function of symbolic powers of a squarefree monomial ideal is eventually constant, one may ask whether the same is true for the Stanley depth. In other words, whether an analogue of Conjecture \ref{limsd} is true, if one replaces the ordinary power with symbolic power. In \cite{s9}, we gave a positive answer to this question. In fact, we have something more. First, we will see in the following theorem that one can compare the Stanley depth of certain symbolic powers of a squarefree monomial ideal.

\begin{thm} [\cite{s9}, Theorem 4.2] \label{sdepsym}
Let $I\subset S$ be a squarefree monomial ideal. Suppose that $m$ and $k$ are positive integers. Then for every integer $j$ with $m-k\leq j\leq m$, we have$${\rm sdepth}(I^{(m)})\geq {\rm sdepth}(I^{(km+j)}) \ \ \ \ \ \ \ \ and \ \ \ \ \ \ \ \ {\rm sdepth}(S/I^{(m)})\geq {\rm sdepth}(S/I^{(km+j)}).$$
\end{thm}

We recall that in the special case of $j=m$, the inequalities of Theorem \ref{sdepsym} were also proved in \cite[Corollary 3.2]{s4}. We also mention that the assertions of Theorem \ref{sdepsym} are true, if one replaces sdepth with depth and this is proved independently by Nguyen and Trung \cite[Theorem 2.7]{nt}, Monta\~no and N\'u\~nez-Betancourt \cite[Theorem 3.4]{mn}, and the author \cite[Theorem 3.3]{s9}.

As an immediate consequences of Theorem \ref{sdepsym}, we obtain the following result.

\begin{cor} [\cite{s9}, Corollary 4.3]
For every squarefree monomial ideal $I\subset S$, we have$${\rm sdepth}(S/I)\geq {\rm sdepth}(S/I^{(2)})\geq {\rm sdepth}(S/I^{(3)})$$and$${\rm sdepth}(I)\geq {\rm sdepth}(I^{(2)})\geq {\rm sdepth}(I^{(3)}).$$
\end{cor}

Assume that $I$ is a squarefree monomial ideal and set$$m:=\min_k{\rm sdepth}(S/I^{(k)}).$$Let $t\geq 1$ be the smallest integer with ${\rm sdepth}(S/I^{(t)})=m$. If $t=1$, then by Theorem \ref{sdepsym}, for every integer $k\geq 1$, we have ${\rm sdepth}(S/I^{(k)})=m$. Now, suppose $t\geq 2$. Again by Theorem \ref{sdepsym}, we have ${\rm sdepth}(S/I^{(t^2-t)})=m$. For every integer $k > t^2-t$, we write $k=st+j$, where $s$ and $j$ are positive integers and $1\leq j\leq t$. As $k > t^2-t$, we conclude that $s\geq t-1$. It then follows from Theorem \ref{sdepsym} that$${\rm sdepth}(S/I^{(k)})={\rm sdepth}(S/I^{(st+j)})\leq {\rm sdepth}(S/I^{(t)})=m.$$By the choice of $m$, we conclude that for every integer $k\geq t^2-t$, the equality ${\rm sdepth}(S/I^{(k)})=m$ holds. Therefore, the sequence $\{{\rm sdepth}(S/I^{(k)})\}_{k=1}^{\infty}$ is convergent and$$\min_k{\rm sdepth}(S/I^{(k)})=m=\lim_{k\rightarrow \infty}{\rm sdepth}(S/I^{(k)}).$$

Similarly, one proves that the sequence $\{{\rm sdepth}(I^{(k)})\}_{k=1}^{\infty}$ is convergent and$$\min_k{\rm sdepth}(I^{(k)})=\lim_{k\rightarrow \infty}{\rm sdepth}(I^{(k)}).$$

Therefore, we have the following result.

\begin{thm} [\cite{s9}, Theorem 4.4] \label{sdsymlim}
For every squarefree monomial ideal $I$, the sequences $\{{\rm sdepth}(S/I^{(k)})\}_{k=1}^{\infty}$ and $\{{\rm sdepth}(I^{(k)})\}_{k=1}^{\infty}$ are convergent. Moreover,$$\min_k{\rm sdepth}(S/I^{(k)})=\lim_{k\rightarrow \infty}{\rm sdepth}(S/I^{(k)}),$$and$$\min_k{\rm sdepth}(I^{(k)})=\lim_{k\rightarrow \infty}{\rm sdepth}(I^{(k)}).$$
\end{thm}

A squarefree monomial ideal $I$ is called {\it normally torsionfree}, if $I^{(k)}=I^k$, for every integer $k\geq 1$. It is immediate from Theorem \ref{sdsymlim} that for any normally torsionfree squarefree monomial ideal $I$, the sequences $\{{\rm sdepth}(S/I^k)\}_{k=1}^{\infty}$ and $\{{\rm sdepth}(I^k)\}_{k=1}^{\infty}$ are convergent. In particular, Conjecture \ref{limsd} is true for normally torsionfree squarefree monomial ideals.

Let $I$ be a squarefree monomial ideal. The smallest integer $t\geq 1$ such that ${\rm depth}(S/I^m)=\lim_{k\rightarrow \infty}{\rm depth}(S/I^k)$ for all $m\geq t$ is called the {\it index of depth stability of powers} of $I$ and is denoted by ${\rm dstab}(I)$. Similarly, one can define the {\it index of depth stability of symbolic powers} by replacing the ordinary powers with symbolic powers. The index of depth stability of symbolic powers is denoted by ${\rm dstab}_s(I)$. By  Theorem \ref{dsymlim}, we have $${\rm dstab}_s(I)\leq n(n+1){\rm bight}(I)^{n/2}.$$According to Theorem \ref{sdsymlim}, one can also define the {\it indices of sdepth stability of symbolic powers}, i.e.,
$${\rm sdstab}_s(I)=\min \big\{t \mid {\rm sdepth}(I^{(m)})=\lim_{k\rightarrow \infty}{\rm sdepth}(I^{(k)}) \ {\rm for \ all} \ m\geq t\big\}$$
$${\rm sdstab}_s(S/I)=\min \big\{t \mid {\rm sdepth}(S/I^{(m)})=\lim_{k\rightarrow \infty}{\rm sdepth}(S/I^{(k)}) \ {\rm for \ all} \ m\geq t\big\}$$We also define the following quantities.
$${\rm sdmin}_s(I)=\min \big\{t \mid {\rm sdepth}(I^{(t)})=\lim_{k\rightarrow \infty}{\rm sdepth}(I^{(k)})\big\}$$
$${\rm sdmin}_s(S/I)=\min \big\{t \mid {\rm sdepth}(S/I^{(t)})=\lim_{k\rightarrow \infty}{\rm sdepth}(S/I^{(k)})\big\}$$

The argument before Theorem \ref{sdsymlim}, also proves the following proposition.

\begin{prop} [\cite{s9}, Corollary 4.5] \label{sstabmin}
For every squarefree monomial ideal $I\subset S$, we have$${\rm sdstab}_s(I)\leq \max\{1, {\rm sdmin}_s(I)^2-{\rm sdmin}_s(I)\}$$and$${\rm sdstab}_s(S/I)\leq \max\{1, {\rm sdmin}_s(S/I)^2-{\rm sdmin}_s(S/I)\}.$$
\end{prop}

As we mentioned above, the assertions of Theorem \ref{sdepsym} are true also for the depth. Thus, a similar argument, as we explained above Theorem \ref{sdsymlim}, implies that the inequalities of Proposition \ref{sstabmin}, remain true, if one replaces Stanley depth with depth. This has been already observed in \cite[Theorem 3.6]{s9}.

Let $I$ be a squarefree monomial ideal. We know from Theorem \ref{sdsymlim} that the sequences $\{{\rm sdepth}(S/I^{(k)})\}_{k=1}^{\infty}$ and $\{{\rm sdepth}(I^{(k)})\}_{k=1}^{\infty}$ are convergent. Now, it is natural to ask the following question.

\begin{ques} \label{qusymlim}
Let $I$ be a squarefree monomial ideal. How can one describe the limits of the sequences $\{{\rm sdepth}(S/I^{(k)})\}_{k=1}^{\infty}$ and $\{{\rm sdepth}(I^{(k)})\}_{k=1}^{\infty}$?
\end{ques}

Question \ref{qusymlim} is widely open. We know the answer only for very special classes of ideals. For example, assume that $I$ is a squarefree complete intersection monomial ideal. It is easy to check that for any integer $k\geq 1$, the equality $I^{(k)}=I^k$ holds. Therefore, using Theorem \ref{compint}, we conclude that$$\lim_{k\rightarrow \infty}{\rm sdepth}(S/I^{(k)})=n-t,$$and$$\lim_{k\rightarrow \infty}{\rm sdepth}(I^{(k)})=n-t+1,$$where $t$ is the number of minimal monomial generators of $I$ (which is also equal to $\ell_s(I)$).

We are also able to compute the limit of the sequence $\{{\rm sdepth}(S/I^{(k)})\}_{k=1}^{\infty}$, where $I$ is the Stanley-Reisner ideal of a matroid. We first recall some basic definitions from the theory of Stanley-Reisner rings.

A {\it simplicial complex} $\Delta$ on the set of vertices $V(\Delta)=[n]:=\{1,
\ldots,n\}$ is a collection of subsets of $[n]$ which is closed under
taking subsets; that is, if $F \in \Delta$ and $F'\subseteq F$, then also
$F'\in\Delta$. Every element $F\in\Delta$ is called a {\it face} of
$\Delta$. The {\it dimension} of a face $F$ is defined to be $|F|-1$. The {\it dimension} of
$\Delta$ which is denoted by $\dim\Delta$, is defined to be $d-1$, where $d
=\max\{|F|\mid F\in\Delta\}$. The {\it Stanley-Reisner ideal} of $\Delta$ is defined as$$I_{\Delta}=\big(\prod_{i\in F}x_i : F\subseteq [n], F\notin \Delta\big)\subseteq S.$$

\begin{dfn}
A simplicial complex $\Delta$ is called {\it matroid} if for every pair of faces $F, G \in \Delta$ with $|F| > |G|$, there is a vertex $x\in F\setminus G$ such that $G\cup\{x\}$ is a face of $\Delta$.
\end{dfn}

As we mentioned above, there are some information about the limit of the Stanley depth function of symbolic powers of Stanley-Reisner ideal of a matroid.

\begin{thm} [\cite{s9}, Theorem 4.7] \label{mat}
Let $\Delta$ be a matroid. Then$$\lim_{k\rightarrow \infty}{\rm sdepth}(S/I_{\Delta}^{(k)})=n-\ell_s(I_{\Delta})=\dim \Delta+1$$and$$\lim_{k\rightarrow \infty}{\rm sdepth}(I_{\Delta}^{(k)})\geq n-\ell_s(I_{\Delta})+1.$$
\end{thm}

\subsection{Cover ideals}

Let $G$ be a graph with vertex set $V(G)=\big\{x_1, \ldots, x_n\big\}$. A subset $C$ of $V(G)$ is called a {\it vertex cover} of $G$ if every edge of $G$ is incident to at least one vertex of $C$. A vertex cover $C$ is called a {\it minimal vertex cover} of $G$ if no proper subset of $C$ is a vertex cover of $G$. The {\it cover ideal} of $G$ is a squarefree monomial ideal of $S$ which is defined as$$J(G)=\big(\prod_{x_i\in C}x_i \mid C \ {\rm is \ a \ minimal \ vertex \ cover \ of} \ G\big).$$It is easy to see that cover ideal is the Alexander dual of edge ideal, i.e.,$$J(G)=I(G)^{\vee}=\bigcap_{\{x_i,x_j\}\in E(G)}(x_1, x_j).$$

Let $I$ be a squarefree monomial ideal. In Question \ref{q1}, we asked whether $I^k$ and $S/I^k$ satisfy Stanley's inequality for every integer $k\gg 0$. One can also ask the similar question for symbolic powers.

\begin{ques} [\cite{s5}, Question 1.2] \label{qssym}
Let $I$ be a monomial ideal. Is it true that $I^{(k)}$ and $S/I^{(k)}$ satisfy Stanley's inequality for every integer $k\gg 0$?
\end{ques}

In this subsection, we investigate the above question for cover ideals, By Theorem \ref{dsymlim}, in order to know whether the high symbolic powers of cover ideals satisfy Stanley's inequality, we need to compute their symbolic analytic spread. This has been done by Constantinescu and Varbaro \cite{cv}. Indeed, they provide a combinatorial description for the symbolic analytic spread of $J(G)$. To state their result, we need to recall some notions from graph theory.

Let $G$ be a graph. A {\it matching} in $G$ is a set of edges such that no two different edges share a common vertex. A subset $W$ of $V(G)$ is called an {\it independent subset} of $G$ if there are no edges among the vertices of $W$. Let $M=\{\{a_i,b_i\}\mid 1\leq i\leq r\}$ be a
nonempty matching of $G$. We say that $M$ is an {\it ordered matching} of
$G$ if the following conditions hold.
\begin{itemize}
\item[(1)] $A:=\{a_1,\ldots,a_r\}$ is an independent subset of vertices of $G$; and

\item[(2)] $\{a_i, b_j\}\in E(G)$ implies that $i\leq j$.
\end{itemize}
The {\it ordered matching number} of $G$, denoted by $\nu_{o}(G)$, is
defined to be $$\nu_{o}(G)=\max\{|M|\mid M\subseteq E(G)\ {\rm is\ an\
ordered\ matching\ of} \ G\}.$$

\begin{thm} [\cite{cv}, Theorem 2.8] \label{sana}
For any graph $G$,$$\ell_s(J(G))=\nu_{o}(G)+1.$$
\end{thm}

As a consequence of Theorems \ref{dsymlim} and \ref{sana}, for any graph $G$ with $n$ vertices we have$$\lim_{k\rightarrow \infty}{\rm depth}(S/J(G)^{(k)})=n-\nu_{o}(G)-1.$$Hoa, Kimura, Terai and Trung \cite{hktt}, determined a linear upper bound for the index of depth stability of symbolic powers of cover ideals. In \cite{s5}, we provided an alternative proof for their result.

\begin{thm} [\cite{hktt}, Theorem 3.4 and \cite{s5}, Theorem 3.1] \label{ldepth}
Let $G$ be a graph with $n$ vertices. Then for every integer $k\geq 2\nu_{o}(G)-1$, we have$${\rm depth}(S/J(G)^{(k)})=n-\nu_{o}(G)-1.$$
\end{thm}

In \cite{s5}, we also proved that high symbolic powers of cover ideals satisfy Stanley's inequality. Indeed, we proved the following result.

\begin{thm} [\cite{s5}, Theorem 3.5 and Corollary 3.6] \label{sdepthsymblim}
Let $G$ be a graph with $n$ vertices. Then for every integer $k\geq 1$, we have
$${\rm sdepth}(J(G)^{(k)})\geq n-\nu_{o}(G) \ \ \  {\rm and} \ \ \  {\rm sdepth}(S/J(G)^{(k)})\geq n-\nu_{o}(G)-1.$$In particular, $J(G)^{(k)}$ and $S/J(G)^{(k)}$ satisfy the Stanley's inequality, for every integer $k\geq 2\nu_{o}(G)-1$.
\end{thm}

The assertions of Theorem \ref{sdepthsymblim} for the special case of bipartite graphs was also proved in \cite{s6}.

Let $G$ be a graph with $n$ vertices. We say $G$ is very well-covered if $n$ is an even integer and moreover, every vertex cover of $G$ has size $n/2$. The graph $G$ is called Cohen-Macaulay if the ring $S/I(G)$ is Cohen-Macaulay. We know from Theorem \ref{sdepthsymblim} that for any graph $G$, the modules $J(G)^{(k)}$ and $S/J(G)^{(k)}$ satisfy the Stanley's inequality, for $k\gg 0$. However, in the case of Cohen-Macaulay very well-covered graphs, we have something more.

\begin{prop} [\cite{s10}, Corollary 3.8] \label{wellsdepth}
Let $G$ be a Cohen-Macaulay very well-covered graph. Then $J(G)^{(k)}$ and $S/J(G)^{(k)}$ satisfy Stanley's inequality, for every integer $k\geq 1$.
\end{prop}

In Question \ref{qusymlim}, we asked about the limit values of the sequences $\{{\rm sdepth}(S/I^{(k)})\}_{k=1}^{\infty}$ and $\{{\rm sdepth}(I^{(k)})\}_{k=1}^{\infty}$, where $I$ is a squarefree monomial ideal. For the case of cover ideals, we pose the following conjecture.

\begin{conj} \label{conjcovsd}
Let $G$ be a graph with $n$ vertices. Then$$\lim_{k\rightarrow \infty}{\rm sdepth}(S/J(G)^{(k)})=n-\nu_{o}(G)-1,$$and$$\lim_{k\rightarrow \infty}{\rm depth}(J(G)^{(k)})=n-\nu_{o}(G).$$
\end{conj}

Let $I$ be a squarefree monomial ideal. According to Theorem \ref{dsymlim}, the sequence $\{{\rm depth}(S/I^{(k)})\}_{k=1}^{\infty}$ is convergent. The situation is even better if $I$ is a cover ideal. In fact, Hoa, Kimura, Terai and Trung \cite[Theorem 3.2]{hktt} proved that the above sequence is non-increasing for cover ideals. In other words, for every graph $G$ and any integer $k\geq 1$, we have$${\rm depth}(S/J(G)^{(k)})\geq {\rm depth}(S/J(G)^{(k+1)}).$$We recall that the above inequality for bipartite graphs was also proved by in \cite[Theorem 3.2]{cpsty}.

We close this article by mentioning that the above inequality is true if one replaces depth with sdepth. In fact, we have the following result.

\begin{thm} [\cite{s5}, Theorem 3.3] \label{nonincsdepthsym}
Let $G$ be a graph. Then for every integer $k\geq 1$, we have

\begin{itemize}
\item[(i)] ${\rm sdepth}(S/J(G)^{(k)})\geq {\rm sdepth}(S/J(G)^{(k+1)})$, and
\item[(ii)] ${\rm sdepth}(J(G)^{(k)})\geq {\rm sdepth}(J(G)^{(k+1)})$.
\end{itemize}
\end{thm}


\section*{Acknowledgment}

The author is grateful to Siamak Yassemi for encouraging him to write this survey article.



\end{document}